\documentclass{article}

\usepackage{amssymb}
\newtheorem{algorithm}{Algorithm}
\newtheorem{definition}{Definition}[section]
\newtheorem{proposition}{Proposition}[section]
\newtheorem{lemma}{Lemma}[section]
\newtheorem{corollary}{Corollary}[section]
\newtheorem{theorem}{Theorem}[section]

\title{\bf On the cardinality of a factor set in the symmetric group}

\author{Krasimir Yordzhev}
\date{}

\begin{document}
\maketitle

\begin{center}
\emph{Faculty of Mathematics and Natural Sciences, South-West University\\
Blagoevgrad,
Bulgaria\\
yordzhev@swu.bg}
\end{center}

\begin{abstract}
Let $n$ be a positive integer,  $\sigma$ be an element of  the symmetric group $\mathcal{S}_n$ and let $\sigma$ be a  cycle of length $n$. The elements $\alpha ,\beta \in \mathcal{S}_n$ are $\sigma$-equivalent, if there are natural numbers $k$ and $l$, such that $\sigma^k \alpha =\beta\sigma^l$, which is the same as the condition to exist natural numbers $k_1$ and $l_1$, such that $\alpha = \sigma^{k_1} \beta\sigma^{l_1}$. In this work we examine some properties of the so defined equivalence relation. We build a finite oriented graph $\Gamma_n$ with the help of which is described an algorithm for solving the combinatorial problem for finding the number of equivalence classes according to this relation.
\end{abstract}



\section{Introduction} \label{intro}

In the present work we denote by $\mathbb{N}$ the set of the nonnegative integers and let $n\in \mathbb{N} $, $n>0$. Let
\begin{equation}\label{zfrak_n}
[n] =\{1,2,\ldots ,n\}.
\end{equation}

We implement the operations $\odot$, $\oplus$ and $\ominus$ which mean multiplication, addition and subtraction by modulus $n$ in $[n]$, i.e. for every $x,y\in[n]$ there are:

\begin{equation}\label{odot}
x\odot y =xy\quad \pmod{n}
\end{equation}
\begin{equation}\label{oplus}
x\oplus y =x+y\quad \pmod{n}
\end{equation}
\begin{equation}\label{ominus}
x\ominus y =x-y\quad \pmod{n}
\end{equation}
$[n]$ together with the operations defined according to (\ref{odot}), (\ref{oplus}) and (\ref{ominus}) is a ring with zero element $n$.

We denote the symmetric group over the set $[n]$ as follows:
\begin{equation}\label{S_n}
\mathcal{S}_n =\left\{ \left. \left(
\begin{array}{cccc}
1    &    2    &  \cdots &    n\\
i_1  &    i_2  &  \cdots &    i_n
\end{array}
\right) \; \right| \; i_s \in [n] ,\; s\in [n],\; i_s \ne i_t \; \mathrm{for} \; s\ne t ,
\right\}
\end{equation}
i.e. the group of all bijective maps $\alpha  :[n] \to [n]$ of the set $[n]$ in itself. We denote by $\alpha (i)\in [n]$ the image of the element $i\in [n]$ under the map $\alpha\in \mathcal{S}_n$.

We assume that the multiplication of  two elements in $\mathcal{S}_n$ is from the left to the right, i.e. if
$$\alpha =\left(
\begin{array}{cccc}
1    &    2    &  \cdots &    n\\
i_1  &    i_2  &  \cdots &    i_n
\end{array}
\right)
\in \mathcal{S}_n ,\quad
\beta =\left(
\begin{array}{cccc}
1    &    2    &  \cdots &    n\\
j_1  &    j_2  &  \cdots &    j_n
\end{array}
\right)
\in \mathcal{S}_n ,$$
then
$$\alpha\beta =\left(
\begin{array}{cccc}
1    &    2    &  \cdots &    n\\
j_{i_1}  &    j_{i_2}  &  \cdots &    j_{i_n}
\end{array}
\right)
\in \mathcal{S}_n .$$

 In other words if  $\alpha ,\beta \in\mathcal{S}_n$ and $\gamma =\alpha\beta$, then for every $i\in [n]$
\begin{equation}\label{alpha*beta}
\gamma (i)=(\alpha\beta)(i)=\beta (\alpha (i)) .
\end{equation}

We denote the identity element of the group $\mathcal{S}_n$ by $\varepsilon$  , i.e.
\begin{equation}\label{varepsilon}
\varepsilon =\left(
\begin{array}{cccc}
1    &    2    &  \cdots &    n\\
1    &    2    &  \cdots &    n
\end{array}
\right)
.
\end{equation}

Let us denote by $ \mathcal{C}_n \subseteq \mathcal{S}_n $ the set of all elements of $\mathcal{S}_n$ which are cycles of length $ n $. Obviously if $\alpha\in \mathcal{C}_n$, then $|\alpha |=n$.

Let
\begin{equation}\label{sigma}
\sigma =\left(
\begin{array}{cccc}
1    &    2    &  \cdots &    n\\
\sigma (1)  &    \sigma (2)  &  \cdots &    \sigma (n)
\end{array}
\right)
\end{equation}

 We suppose that we have chosen and fixed $\sigma\in \mathcal{C}_n$ described using a formula
 (\ref{sigma}) and to the end of the article when we write $\sigma$ we will understand the concrete chosen in this moment element of $\mathcal{C}_n$.

\begin{definition}
Let $\alpha ,\beta ,\sigma \in \mathcal{S}_n$. We say that $\alpha$ and $\beta$ are $\sigma$-equivalent (or just equivalent, if  $\sigma$ is understood) and we write $\alpha \stackrel{\sigma}{\sim} \beta$ (or just $\alpha \sim\beta$ if $\sigma$ is understood), if there are $k,l\in \mathbb{N}$, such that $\sigma^k \alpha =\beta\sigma^l$, which is the same as the condition  to exist $k_1 ,l_1 \in \mathbb{N}$, such that $\alpha = \sigma^{k_1} \beta\sigma^{l_1}$
\end{definition}

Obviously, this relation is an equivalence relation.

If  $\alpha\in \mathcal{S}_n$, then  we denote by $K_\alpha$ the set of elements of $\mathcal{S}_n$ to which $\alpha$ is related by the equivalence relation $\stackrel{\sigma}{\sim}$, i.e.

\begin{equation}\label{K_alpha}
K_\alpha =\{ \beta \in \mathcal{S}_n \;|\; \exists k,l\in \mathbb{N} :\; \sigma^k \alpha=\beta \sigma^l \}.
\end{equation}

The equivalence classes of $\mathcal{S}_n$ by the equivalence  relation
$\stackrel{\sigma}{\sim}$ are a particular kind of  {\it double coset} \cite{bogopolski,curtis,vos}. They make use of substitutions group theory  and linear representation of finite groups theory \cite{curtis,vos}.

If $$\sigma =
\left(
\begin{array}{cccccc}
1    &    2    &  3  &  \cdots &    n-1   &   n\\
2    &    3    &  4  &  \cdots &    n     &   1
\end{array}
\right) ,$$ or $$\sigma =
\left(
\begin{array}{cccccc}
1    &    2    &  3  &  \cdots &    n-1   &   n\\
n    &    1    &  2  &  \cdots &    n-2   &   1
\end{array}
\right) ,$$ then $\sigma$-equivalence take place in the mathematical modeling in the textile industry \cite{borzunov,bansko,umb2010}.

If $s$ and $t$ are positive integers, then we denote by $\mathrm{GCD}(s,t)$ the greatest common divisor of $s$ and $t$.

With  $\varphi (n)$ we denote the Euler totient function, i.e. the number of elements of $[n]$  that are relatively prime to $n$ \cite{andreescu,mirchevNumber}. By definition  $\varphi (1)=1$.

Let $\sigma\in \mathcal{C}_n$. The aim of the present work is to describe the effective algorithm, which solves the combinatoric problem of finding the number of equivalence classes according to $\sigma$-equivalence when $n$ is given, or in other words to find the number of the elements in {\it the factor-set}
\begin{equation}\label{Qn}
Q_n = {\mathcal{S}_n }_{/_{\stackrel{\sigma}{\sim}}}.
\end{equation}

As we will see in the presentation the number is one and the same for any $\sigma \in \mathcal{C}_n$ and this number depends only of $n$.

In our work we will use some concepts from graph theory \cite{diestel,harary,mirchevGraphs}.

\section{Some properties of the $\sigma$-equivalence}\label{sec2}

 Let $\sigma\in\mathcal{C}_n$. For our aim it is necessary to find the number of all solutions of the equation
\begin{equation}\label{eq}
\sigma^k \xi =\xi \sigma^l ,
\end{equation}
where
\begin{equation}\label{x}
\xi =\left(
\begin{array}{cccc}
1    &    2    &  \cdots &    n\\
x_1  &    x_2  &  \cdots &    x_n
\end{array}
\right)
\in \mathcal{S}_n .
\end{equation}

\begin{proposition} \label{prop1}
Let $n\in \mathbb{N}$, $\sigma ,\xi \in \mathcal{S}_n$, defined according to (\ref{sigma}) and (\ref{x}). Then the following equalities are true:
\begin{equation}\label{sigma*x}
\sigma \xi =\left(
\begin{array}{cccc}
1    &    2    &  \cdots &    n\\
x_{\sigma (1)}  &    x_{\sigma (2)}  &  \cdots &    x_{\sigma (n)}
\end{array}
\right)
\end{equation}
\begin{equation}\label{x*sigma}
\xi \sigma =\left(
\begin{array}{cccc}
1    &    2    &  \cdots &    n\\
\sigma (x_1 )  &    \sigma (x_2 )  &  \cdots &    \sigma (x_n )
\end{array}
\right)
\end{equation}
\begin{equation}\label{sigma^k*x}
\sigma^k \xi =\left(
\begin{array}{cccc}
1    &    2    &  \cdots &    n\\
x_{\sigma^k (1)}  &    x_{\sigma^k (2)}  &  \cdots &    x_{\sigma^k (n)}
\end{array}
\right)
\end{equation}
\begin{equation}\label{x*sigma^l}
\xi \sigma^l =\left(
\begin{array}{cccc}
1    &    2    &  \cdots &    n\\
\sigma^l (x_1 )  &    \sigma^l (x_2 )  &  \cdots &    \sigma^l (x_n )
\end{array}
\right)
\end{equation}
\end{proposition}

The proof follows from the definition of the operation multiplication in the symmetric group $\mathcal{S}_n$ and according formulas (\ref{alpha*beta}), (\ref{sigma}) and (\ref{x}).

\hfill $\Box$

\begin{lemma} \label{lm1}
Let $n\in \mathbb{N}$, $\xi ,\sigma \in \mathcal{S}_n$. Then for every number $s\in \mathbb{N}$ from the equation $\sigma^k \xi =\xi \sigma^l$ follows the equation $\sigma^{sk} \xi =\xi \sigma^{sl}$
\end{lemma}

Proof. Induction by $s$. When $s=0$ and $s=1$ the statement is obviously true. We assume that $\sigma^k \xi =\xi \sigma^l$ and $\sigma^{(s-1)k} \xi =\xi \sigma^{(s-1)l}$, $s\in \mathbb{N}$, $s>1$ are true. Then, $\sigma^{sk} \xi =\sigma^k \sigma^{(s-1)k} \xi =\sigma^k \xi \sigma^{(s-1)l} =\xi \sigma^l \sigma^{(s-1)l} =\xi \sigma^{sl} $

\hfill $\Box$

\begin{corollary} \label{crl1}
The set of all the solutions of the equation $\sigma^{sk} \xi =\xi \sigma^{sl} ,$ $\sigma ,\xi  \in \mathcal{S}_n $ contains all the solutions of the equation $\sigma^k \xi  =\xi \sigma^l $ too.

\hfill $\Box$
\end{corollary}

\begin{lemma} \label{lm3}
Let $n\in \mathbb{N}$, $\sigma \in \mathcal{S}_n$ and let $|\sigma |=n$. Then if $K_\xi  \in {\mathcal{S}_n }_{/_{\stackrel{\sigma}{\sim}}}$ is an arbitrary equivalence class according to the relation $\stackrel{\sigma}{\sim}$ and  $\xi\in\mathcal{S}_n$ is its representative, then the following two cases are possible:

$i)$ $| K_\xi  | =n^2 $

$ii)$ $\xi \in S_n$ is a solution of a equation of the following kind:
\begin{equation}\label{4}
\sigma^k \xi =\xi \sigma^l
\end{equation}
and the following conditions come true:
\begin{equation}\label{5}
1\le k\le l <n
\end{equation}
\begin{equation}\label{6}
k|n
\end{equation}
\begin{equation}\label{7}
|K_\xi  |=kn
\end{equation}
\begin{equation}\label{k|l}
k|l
\end{equation}
\begin{equation}\label{8}
\exists s\in \mathbb{N}:\;1\le s<n,\;  \mathrm{GCD} (s,n)=1, \; l=s\odot k
\end{equation}
\end{lemma}

Proof. We can get all elements of $K_\xi $ by building the table:
$$\begin{array}{rrrrr}
\xi  & \xi \sigma & \xi \sigma^2 & \cdots & \xi \sigma^{n-1} \\
\sigma \xi  & \sigma \xi  \sigma & \sigma \xi  \sigma^2 & \cdots & \sigma \xi  \sigma^{n-1} \\
\sigma^2 \xi  & \sigma^2 \xi \sigma & \sigma^2 \xi \sigma^2 & \cdots & \sigma^2 \xi \sigma^{n-1} \\
\cdots \\
\sigma^{n-1} \xi  & \sigma^{n-1} \xi  \sigma & \sigma^{n-1} \xi  \sigma^2 & \cdots & \sigma^{n-1} \xi \sigma^{n-1} \\
\end{array} $$

Since $|\sigma |=n$, then it is easy to see that all elements in a given row or column of the obtained table are different.

For every $\xi \in S_n$ the following two cases are possible:

$i)$ There are not natural numbers $k,l$, such that $\sigma^k \xi =\xi \sigma^l$. Then obviously, all elements in the table are different and it follows that $| K_\xi  | =n^2 $.

$ii)$ There are natural numbers $k,l$ such that $\sigma^k \xi =\xi \sigma^l $
and as (\ref{4}) is equal to $\sigma^{n-k} \xi =\xi \sigma^{n-l}$, then without restriction of the generality we can consider that $1\le k\le l <n $. Let $k$ is the minimal natural number, for which (\ref{4}) and (\ref{5}) are true when $\xi \in S_n$ is given. In this case we can see that $k$ divides $n $ and $|K_\xi  |=kn $. We proved (\ref{5}), (\ref{6}) and (\ref{7}).

According to (\ref{6}) there is $p\in N$, such that $n=pk$. From lemma \ref{lm1} and from $\sigma^k \xi =\xi \sigma^l$ it follows that $\sigma^{pk} \xi =\xi \sigma^{pl}$. But $\sigma^{pk} =\sigma^n =\varepsilon \Longrightarrow \sigma^{pl} =\varepsilon$. Therefore, $n=pk$ divides $pl$, from where $k$ divides $l$, i.e. there is $q\in N$, such that $l=qk$. We proved (\ref{k|l}).

If $\mathrm{GCD}(q,n)=1$, then $s=q$ and in this case (\ref{8}) is true.

Let $\mathrm{GCD}(q,n)>1$, where $\displaystyle q=\frac{l}{k}$ and let $\displaystyle p=\frac{n}{k}$. If we assume that $\mathrm{GCD}(p,q)=t>1$, then $p=p_1 t,$ $q=q_1 t$. Obviously, $k$ can be represented in the kind $k=uv$, where $u=u_1^{a_1} u_2^{a_2} \ldots u_w^{a_w} ,$ $u_i$-prime number, $u_i |q$, $i=1,2,\ldots ,w$ and $\mathrm{GCD}(v,q)=1$. Then we set $s=n-vp+q$. It is easy to see that $1\le s<n$ and $ks=kn-kpv +kq\equiv kq=l \pmod{n}$. Therefore in this case $l=s\odot k$. We have to prove that in that case $\mathrm{GCD}(s,n)=1$ too.

Let $r\in N,$ $r>1$ is a prime number and let $r$ divides $n$.

If $r$ divides $p$, then $r$ does not divides $q$ $\Longrightarrow$ $r$ does not divides $s =n-vp+q$.

If $r$ does not divides $p$, then $r$ divides $k$ and there are two cases

a) $r$ divides $u$ and $r$ does not divides $v$. Then $r$ divides $q$, but $r$ does not divides $vp$ $\Longrightarrow$ $r$ does not divides $s=n-vp+q$.

b) $r$ does not divides $u$ and $r$ divides $v$. Then $r$ does not divides $q$ $\Longrightarrow$ $r$ does not divides $s=n-vp+q$.

Hence, $\mathrm{GCD}(s,n)=1$. We proved that (\ref{8}) is true.

\hfill $\Box$

\begin{lemma}\label{lm4}
Let $n\in \mathbb{N}$, $\sigma\in\mathcal{S}_n$, $|\sigma |=n$. Then if $\xi ,\eta\in S_n$ are such that $\sigma^k \xi =\xi \sigma^l ,$ $\sigma^k \eta=\eta\sigma^r ,$ $l\not\equiv r \pmod{n}$, then $K_\xi \ne K_\eta $
\end{lemma}

Proof. We assume the opposite thing, i.e. having in mind lemma \ref{lm3} we assume that there are $\xi ,\eta\in S_n$ and natural numbers $k,l,r$, such that $1\le k,l,r<n$ and for which $\sigma^k \xi=\xi\sigma^l ,$ $ \sigma^k \eta=\eta\sigma^r$ and $K_\xi =K_\eta $ are true. Then there are natural numbers $p$ and $q$, such that $\eta=\sigma^p \xi \sigma^q$. Replacing we get $\sigma^{k+p} \xi\sigma^q =\sigma^p \xi\sigma^{q+r} \Longrightarrow $ $\sigma^k \xi=\xi\sigma^r \Longrightarrow $ $\xi\sigma^l =\xi\sigma^r \Longrightarrow $ $\sigma^l =\sigma^r$ and since $1\le l,r< n$, and $|\sigma | =n$, then $l=r$, which is in conflict with the condition that $l\not\equiv r \pmod{n}$.

\hfill $\Box$
\begin{lemma} \label{lm2}
Let $n\in \mathbb{N}$, $\sigma \in \mathcal{C}_n$ and $\xi \in \mathcal{S}_n$,  defined using respectively formulas
(\ref{sigma}) and (\ref{x}). Then the equation
\begin{equation}\label{sx=xs^l}
\sigma \xi =\xi \sigma^l
\end{equation}
has a solution in $\mathcal{S}_n$, if and only if $\mathrm{GCD}(l,n)=1 $.
\end{lemma}

Proof. Let $\xi \in \mathcal{S}_n$ is a solution of the equation $\sigma \xi =\xi\sigma^l $. We assume that there are whole positive numbers $l_1 ,n_1$ and $p$, such that $1\le n_1 <n$, $1<p\le n$,  $l=l_1 p$ and $n=n_1 p$, i.e. $\mathrm{GCD}(l,n)=p>1$. Then from lemma \ref{lm1} it follows that $\sigma^{n_1} \xi =\xi \sigma^{l_1 pn_1 } =\xi \sigma^{l_1 n } =\xi \varepsilon =\xi $, from where  $\sigma^{n_1} =\varepsilon$, which is not possible, as $n_1 <n$, and $|\sigma| = n$. So we proved the necessity.

Sufficiency. Let $a$ is an arbitrary element of $[n]$. Having in mind (\ref{sigma*x}) and (\ref{x*sigma^l}) we will look for solution of the system of $n$ equation with $n$ unknowns $x_1 ,x_2 ,\ldots ,x_n$
\begin{equation}\label{syst}
\| \; x_{\sigma (i) } =\sigma^l (x_i ) , \quad i=1,2,\ldots ,n
\end{equation}
as the beginning condition is
\begin{equation}\label{x1}
x_1 =a.
\end{equation}

We will prove that there are $x_1 ,x_2 , \ldots ,x_n \in [n]$ satisfying system (\ref{syst}) as the beginning condition is (\ref{x1}), then for every $k\in \mathbb{N}$
\begin{equation}\label{recurs}
x_{\sigma^k (1)} =\sigma^{kl} (a) .
\end{equation} is true.
When $k=0$ (\ref{recurs}) follows from the beginning condition (\ref{x1}), and when $k=1$ it follows from (\ref{syst}) as $i=1$ and from (\ref{x1}). We assume that (\ref{recurs}) is true for some $k\in \mathbb{N}$. Then we set $i=\sigma^k (1)\in [n]$ and as we reckon with (\ref{syst}) we consecutively get: $x_{\sigma^{k+1} (1)} =x_{\sigma (\sigma^k (1))} =x_{\sigma (i)} = \sigma^l (x_i )= \sigma^l (x_{\sigma^k (1)} )=\sigma^l (\sigma^{kl} (a))=\sigma^{(k+1)l} (a)$. Hence, the equation (\ref{recurs})  is true for every $k\in \mathbb{N}$.

Since $\sigma\in\mathcal{C}_n$, then the set $\{ x_1 ,x_{\sigma (1)} ,x_{\sigma^2 (1)} ,\ldots ,x_{\sigma^{n-1} (1)} \}$ is identical with the set  $X=\{ x_1 ,x_2 ,\ldots ,x_n \}$. Hence, using the recursive dependence (\ref{recurs}) as the beginning condition is (\ref{x1}) when $k=1,2,\ldots ,n-1$ we definitely get the solution of the system (\ref{syst}) as the condition is (\ref{x1}).

We have to prove that if $\mathrm{GCD} (n,l)=1$, then $$\xi  =\left(
\begin{array}{cccc}
1    &    2    &  \cdots &    n\\
x_1  &    x_2  &  \cdots &    x_n
\end{array}
\right) \in \mathcal{S}_n ,$$
 i.e. $x_s \ne x_t$ when $s\ne t$. We assume the opposite, i.e. we assume that there are $s,t\in [n]$, $s>t$ such that $x_{\sigma^s (1)} =x_{\sigma^t (1)}$. From lemma \ref{lm1} it follows the equation $\sigma^{sl} (a) =\sigma^{tl} (a) $, from where it follows that $\sigma^{(s-t)l} (a)=\varepsilon (a)=a$. But $a\in [n]$ is randomly chosen and the last equation is true for every $a\in [n]$. This is possible, if and only if, $(s-t)l \equiv n \pmod{n}$. And since $\mathrm{GCD} (n,l)=1$, then it follows that $(s-t)\equiv n\pmod{n}$ is true, which is impossible, because $1\le s,t\le n$, from where $0<s-t<n$.
 Therefore, $\xi  =\left(
\begin{array}{cccc}
1    &    2    &  \cdots &    n\\
x_1  &    x_2  &  \cdots &    x_n
\end{array}
\right) \in \mathcal{S}_n$.

\hfill $\Box$

 Lemma \ref{lm2} is useful, because in the so proposed proof is described an algorithm for finding all solutions in $\mathcal{S}_n$ of the equation  (\ref{sx=xs^l}) in case that $\mathrm{GCD} (n,l)=1$.

\section{A graph used for counting the number of all elements of the factor-set ${\mathcal{S}_n }_{/_{\stackrel{\sigma}{\sim}}}$ }\label{sec3}

Using the base of the propositions, described in section \ref{sec2} we build finite oriented graph $\Gamma_n =\left( V_n ,\Psi\right)$ consisting of set of vertices $V_n \subset [n] \times [n] = \{ \langle k,l\rangle \; |\; 1\le
k,l\le n\}$ and set of arcs $\Psi \subset V_n \times V_n$. This graph helps us to find all equivalence classes according to the equivalence relation $" \stackrel{\sigma}{\sim} " $, where $n\in \mathbb{N}$, $\sigma \in \mathcal{C}_n$.  The following  algorithm helps us to build $\Gamma_n$:

\begin{algorithm} \label{alg1}
Building oriented graph $\Gamma_n$:
\begin{enumerate}
\item For every natural number $l$, such that $1\le l<n$ and $\mathrm{GCD}(l,n)=1$ we build vertex $\langle 1,l\rangle \in V_n$.
\item If the vertex $\langle k,l\rangle \in V_n $ is already built, and
$p$ is a prime divisor of $n$, such that $kp\; |\; n$, then we build vertex $\langle kp, l\odot p\rangle \in V_n$ (if it is not built yet) and we build arc  $\left( \langle k,l\rangle ,\langle
kp, l\odot p\rangle \right)$ with the beginning $\langle k,l\rangle$ and end  $\langle{kp, l\odot p}\rangle$.

\item We repeat point 2 until the possibilities are exhausted.
\end{enumerate}
\end{algorithm}

Point 1 of algorithm \ref{alg1} is made according to the propositions in lemma \ref{lm3} and lemma \ref{lm2}. Building the vertices of $V_n$ in point 2 of algorithm \ref{alg1} is according the propositions in lemma \ref{lm3}, lemma \ref{lm1} and corollary \ref{crl1}, and the arcs are built according lemma \ref{lm1} and corollary \ref{crl1}. There are not other vertices in $\Gamma_n $, besides these received with the help of algorithm \ref{alg1}, according to lemma \ref{lm4}.

Studying algorithm \ref{alg1} it is easy to convince that for every $n\in \mathbb{N}$ the graph $\Gamma_n$ is connected, there is a vertex $\langle n,n\rangle \in V_n$, which is accessible  from every one of the other vertices of the graph, and this vertex is not a beginning of any arc. There is not arc ending in the vertices built in point 1 of algorithm \ref{alg1} of the kind $\langle 1,l\rangle$.

We naturally define  partial order $\prec$ in the set of vertices $V_n$, namely $\langle s,t\rangle \prec \langle k,l\rangle$, if and only if, there is a directed path in $\Gamma_n$, where the first vertex is $\langle s,t\rangle$ and the last is $\langle k,l\rangle$. Obviously, $V_n $ together with the so initiated partial order is semilattice with maximal element $\langle n,n\rangle \in V_n$ and $\varphi (n)$ in number minimal elements, each of kind $\langle 1, l\rangle$, where  $\varphi (n)$ is the Euler totient function.

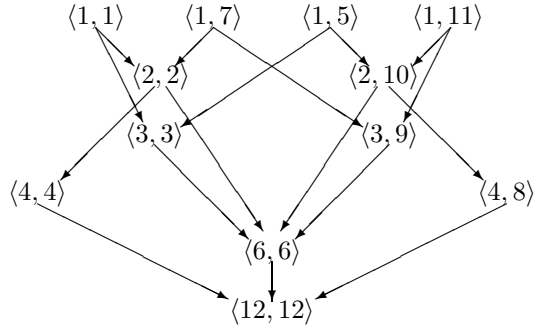
\begin{figure}[h]
\unitlength=3.9mm\linethickness{0.4pt}
\begin{center}
\begin{picture}(16.00,12.50)
\put(2.00,10.00){\makebox(0,0)[cc]{$\langle{1,1}\rangle$}}
\put(6,10){\makebox(0,0)[cc]{$\langle{1,7}\rangle$}}
\put(10,10){\makebox(0,0)[cc]{$\langle{1,5}\rangle$}}
\put(14,10){\makebox(0,0)[cc]{$\langle{1,11}\rangle$}}
\put(4.2,8){\makebox(0,0)[cc]{$\langle{2,2}\rangle$}}
\put(11.8,8){\makebox(0,0)[cc]{$\langle{2,10}\rangle$}}
\put(2,9.66){\vector(1,-1){1.3}} \put(6,9.66){\vector(-1,-1){1.3}}
\put(10,9.66){\vector(1,-1){1.3}}
\put(14.1,9.66){\vector(-1,-1){1.3}}
\put(4,6){\makebox(0,0)[cc]{$\langle{3,3}\rangle$}}
\put(12,6){\makebox(0,0)[cc]{$\langle{3,9}\rangle$}}
\put(4,5.66){\vector(1,-1){3.2}}
\put(12,5.66){\vector(-1,-1){3.2}}
\put(2,9.66){\vector(1,-2){1.6}}
\put(14.1,9.66){\vector(-1,-2){1.6}}
\put(6,9.66){\vector(3,-2){5.1}}
\put(10,9.66){\vector(-3,-2){5.1}}
\put(0,4){\makebox(0,0)[cc]{$\langle{4,4}\rangle$}}
\put(16,4){\makebox(0,0)[cc]{$\langle{4,8}\rangle$}}
\put(4,7.66){\vector(-1,-1){3.2}}
\put(12,7.66){\vector(1,-1){3.2}}
\put(4.4,7.66){\vector(2,-3){3.3}}
\put(11.6,7.66){\vector(-2,-3){3.3}}
\put(8,2){\makebox(0,0)[cc]{$\langle{6,6}\rangle$}}
\put(0,3.66){\vector(2,-1){6.5}}
\put(16,3.66){\vector(-2,-1){6.5}}
\put(8,0){\makebox(0,0)[cc]{$\langle{12,12}\rangle$}}
\put(8,1.66){\vector(0,-1){1.33}}
\end{picture}
\caption{The graph $\Gamma_{12}$}\label{f3}
\end{center}
\end{figure}

For example on figure  \ref{f3} are shown the graph $\Gamma_{12}$.

\begin{theorem}\label{pnk}
Let $n\in \mathbb{N}$, $\sigma \in \mathcal{C}_n$ and let the conditions (\ref{5}) $\div$ (\ref{8}) are satisfied. Then the number of all solutions of the equation (\ref{4})  is equal to
\begin{equation}\label{fpnk}
{\rm p}(n,k) =k! \left( \frac{n}{k} \right)^k .
\end{equation}
\end{theorem}

Proof. Consecutively we get the sets $M=\left\{ m_1 ,m_2 ,\ldots ,m_k \right\} $ and $A_1 ,A_2 ,\ldots ,A_k$, where $M\subset [n] $, $A_i \subset [n]$, $i=1,2,\ldots ,k$, using the following algorithm:

When $i=1$ we set $m_1 =1$ and
$$
\begin{array}{rcl}
A_1 & = & \left\{ m_1 , \sigma^k (m_1 ),\sigma^{2k} (m_1 ),\ldots ,\sigma^{n-k} (m_1 )\right\} = \\
    & = & \left\{ \sigma^{tk} (m_1 )\; | \; t=0,1,\ldots ,\frac{n}{k} -1 \right\}
\end{array}
$$

Let the numbers $m_1 ,m_2 , \ldots ,m_{i-1} \in M$ are received, and let the sets $A_1 ,A_2 , \ldots A_{i-1}$, $1<i\le k$ are received too. Then we set:
\begin{equation}\label{mi}
m_i =\min \left( [n] \setminus \left( \bigcup_{j=1}^{i-1} A_j \right) \right)
\end{equation}

\begin{equation}\label{Ai}
\begin{array}{rcl}
A_i & = & \left\{ m_i , \sigma^k (m_i ),\sigma^{2k} (m_i ),\ldots ,\sigma^{n-k} (m_i )\right\} = \\
    & = &\left\{ \sigma^{tk} (m_i )\left|  t=0,1,\ldots ,\frac{n}{k} -1 \right. \right\}
\end{array}
\end{equation}

 When $\displaystyle t=\frac{n}{k}$ it is satisfied $\sigma^{tk} (m_i )=\sigma^n (m_i )=m_i$, $i=1,2,\ldots ,k$, i.e. $\sigma^{tk} (m_i )\in A_i $ for every $t\in \mathbb{N} $ and $A_i$ is built to contain just numbers of the kind $\sigma^{tk} (m_i )$, $t\in \mathbb{N}$. It follows that $$\left| A_i \right| =\frac{n}{k} \quad i=1,2,\ldots ,k.$$

Since $\sigma \in \mathcal{C}_n$ and $\sigma^r (a) \ne \sigma^r (b)$ as $a\ne b$, where $1\le a,b,r\le n$, then the received family of subsets $A_1 ,A_2 ,\ldots ,A_k$ is a partition of $[n]$, i.e. $A_i \cap A_j =\emptyset$ as $i\ne j$ and $\displaystyle \bigcup_{i=1}^k A_i =[n]$.

Let $a\in [n]$ and let more concrete $a\in A_j$ for some $j\in \{ 1,2,\ldots ,k\}$. We examine the set $$B=\left\{ \sigma^{tl} (a)\; |\; t=0,1,\ldots ,\frac{n}{k} -1 \right\}$$
We will prove that $B=A_j$. And really, according to (\ref{8}) there is natural number $s$, such that $1\le s< n$ and $l=s\odot k$. Then for every $\displaystyle t\in \left\{ 0,1,\ldots , \frac{n}{k}-1 \right\}$ there is $\sigma^{tl} (a)=\sigma^{t\odot l} (a) = \sigma^{t\odot s\odot l } (a)=\sigma^{t_1 k} (a) \in A_j$ for some  $\displaystyle t_1 \in \left\{ 0,1,\ldots ,\frac{n}{k} -1\right\}$. Therefore, $B\subseteq A_j$. We have to prove that $\displaystyle |B|= \frac{n}{k}$, from where it follows that $B=A_j$. We assume the opposite thing, i.e. we assume that there are natural numbers $t_1 ,t_2$, such that $\displaystyle 0\le t_1 <t_2 < \frac{n}{k}$ for which $\sigma^{t_2 l} (a) =\sigma^{t_1 l} (a)$ is true. It follows that $\sigma^{(t_2 -t_1)l} (a)=\varepsilon (a) =a$. But $\sigma$ is a cycle with length $n$ as it is given and therefore the last equation is possible, if and only if, $(t_2 -t_1 )l \equiv 0\pmod{n}$. But according to (\ref{8}) $l=s\odot k$, where $ \mathrm{GCD} (s,n)=1$. Therefore the last comparison is equal to the comparison $(t_2 -t_1 )k\equiv 0 \pmod{n}$. But  $\displaystyle 0\le t_1 <t_2 <\frac{n}{k}$, from where $\displaystyle 0< t_2 -t_1 <\frac{n}{k}$ and because $t_1 \ne t_1$, then $k=0$, which is impossible
. Therefore $B=A_j$.

Having in mind proposition \ref{prop1} it follows that $$\displaystyle \xi =\left(
\begin{array}{cccc}
1    &    2    &  \cdots &    n\\
x_1  &    x_2  &  \cdots &    x_n
\end{array}
\right)
\in \mathcal{S}_n$$
is a solution of the equation $\sigma^k \xi =\xi\sigma^l$, if and only if, $x_1 ,x_2 ,\ldots ,x_n$ are solutions of the system of $n$ equations with $n$ unknowns
\begin{equation}\label{c-ma}
\| x_{\sigma^k (r)} =\sigma^l (x_r ),\quad r=1,2,\ldots ,n
\end{equation}
where $x_i \in [n]$, $i=1,2,\ldots ,n$, $x_i \ne x_j$ as $i\ne j$.

We part the set of unknowns $X=\{ x_1 ,x_2 ,\ldots , x_n \}$ into the subsets $X_1 ,X_2 ,\ldots ,X_k$, where
$$\begin{array}{rcl}
X_i & = &\left\{ x_a \; |\; a\in A_i \right\} = \\
    & = &\left\{ x_{m_i} ,x_{\sigma^k (m_i )} ,x_{\sigma^{2k} (m_i )} ,\ldots ,x_{\sigma^{n-k} (m_i )} \right\} ,\quad  i=1,2,\ldots ,k
\end{array}$$

Obviously, for every $i\in \{ 1,2,\ldots ,k\}$ $\displaystyle |X_i |=\frac{n}{k}$, $X_i \cap X_j =\emptyset$ as $i\ne j$ and $\displaystyle \bigcup_{i=1}^k X_i =X$.

Let $a\in A_j$ for some $j\in \{ 1,2,\ldots ,k\}$ and let $i\in \{ 1,2,\ldots ,k\}$. We set
$$x_{m_i} =a$$
Then according to (\ref{c-ma}) we consecutively get:
$$x_{\sigma^k (m_i )} =\sigma^l (x_{m_i } ) =\sigma^l (a)$$
$$x_{\sigma^{2k} (m_i )} =x_{\sigma^k (\sigma^k (m_i ))} =\sigma^l (\sigma^k (m_i )) =\sigma^l (\sigma^l (a)) =\sigma^{l\odot 2} (a)$$

Continue do this by induction we get that
$$x_{\sigma^{tk} (m_i )} =\sigma^{l\odot t } (a) \quad \forall t\in \left\{ 0,1,2,\ldots ,\frac{n}{k} -1\right\}$$

According to the above proved, we can get concrete values for every unknown of the set $X_i$, as all values are different to each other and they are among the elements of the set $A_j \ni a$. Hence, the problem for finding all solutions of the system (\ref{c-ma}) is put to the problem for choosing the elements $a_1 \in A_{i_1}$, $a_2 \in A_{i_2}$, $\ldots$, $a_k \in A_{i_k}$, where $i_r \ne i_s$ as $r\ne s$. This can be made by $\displaystyle k! \left( \frac{n}{k} \right)^k$ ways, which proves the theorem.

\hfill $\Box$

\begin{corollary}\label{cpnk}
Let $n\in \mathbb{N}$, $\sigma \in \mathcal{C}_n$ and let the conditions (\ref{5}) $\div$ (\ref{8}) are satisfied. Then the number of all solutions of the equation (\ref{4}) depends only of $n$ and $k$ and it does not depends of $l$.

\hfill $\Box$
\end{corollary}

Having in mind lemma \ref{lm2} and theorem \ref{pnk} we get:

\begin{corollary}\label{kk=1}
If $\mathrm{GCD}(l,n)=1 $, then the equation (\ref{sx=xs^l}) has unique solution in $\mathcal{S}_n$ to within $\sigma$-equivalence relation.

\hfill $\Box$
\end{corollary}

When we proved the theorem \ref{pnk} we described an algorithm for finding all solutions of the equation (\ref{4}) when (\ref{5}) $\div$ (\ref{8}) are true.

From the definition of $\Gamma_n =\left( V_n ,
\Psi\right)$ (algorithm \ref{alg1}), lemma  \ref{lm3}, lemma \ref{lm4} and lemma \ref{lm2} it follows that there is surjective map $f:V_n \to {\mathcal{S}_n }_{/_{\stackrel{\sigma}{\sim}}}$ from the set of vertices  $V_n$  on the set ${\mathcal{S}_n }_{/_{\stackrel{\sigma}{\sim}}}$, defined as follows:

If $\langle k, l\rangle \in V_n$, where $k< n$, then

\begin{equation}\label{phik}
f (\langle k,l\rangle )=\left\{ K_\xi \in {\mathcal{S}_n }_{/_{\stackrel{\sigma}{\sim}}} \; \left| \;  \xi \in S_n ,\; \sigma^k \xi =\xi \sigma^l ,\;  \sigma^{k_1} \xi \ne \xi \sigma^{l_1} \; \forall k_1 <k \right. \right\}
\end{equation}
and
\begin{equation}\label{phin}
f (\langle n,n\rangle )=\left\{ K_\xi \in {\mathcal{S}_n }_{/_{\stackrel{\sigma}{\sim}}} \; \left|\; |K_\xi |=n^2 \right. \right\}
\end{equation}

\begin{theorem}\label{t5}
Let $n\in \mathbb{N}$, $\sigma\in \mathcal{C}_n$ and let
\begin{equation}\label{hnk}
\begin{array}{rcl}
 {\rm h}(n,k) & = & \left| f (\langle k,l\rangle ) \right| = \\
 & = & \left| \left\{ K_\xi \in {\mathcal{S}_n }_{/_{\stackrel{\sigma}{\sim}}} \; |\; \exists \langle k,l\rangle \in V_n \; : \;  K_\xi \in f (\langle k,l\rangle ) \right\} \right| ,
 \end{array}
\end{equation}
where $f : V_n \to {\mathcal{S}_n }_{/_{\stackrel{\sigma}{\sim}}}$ is the map defined with the help of (\ref{phik}) and (\ref{phin}). Then for every couple of natural numbers $k$ and $l$, such that  $\langle k,l \rangle$ is a vertex of the graph $\Gamma_n$:
\begin{equation}\label{fhnk}
{\rm h}(n,k) =\frac{1}{k} \left( (k-1)! \left( \frac{n}{k} \right)^{k-1} - \sum_{\langle r,t\rangle \prec \langle k,l\rangle}  r\, {\rm h}(n,r) \right)
\end{equation}
is true.
\end{theorem}

Proof. First we have to emphasize that from theorem \ref{pnk} it follows that $h(n,k)$ does not depend of the parameter $l$ and it is one and the same for every vertex $\langle k, l\rangle \in V_n$.

Let $\langle k, l\rangle \in V_n$ and let $1\le k< n$. Having in mind corollary \ref{crl1} and lemma \ref{lm3} we get:
$$kn\, {\rm h} (n,k) +\sum_{\langle r,t\rangle \prec \langle k,l\rangle} rn\, {\rm h} (n,r) = {\rm p}(n,k) ,$$
where ${\rm p}(n,k) $ is the function, defined using the formulation of the theorem \ref{pnk}, i.e. ${\rm p}(n,k)$ is equal to the number of all solutions of the equation $\sigma^k \xi =\xi \sigma^l$ as the conditions (\ref{5})$\div$(\ref{8}) are true.  Then applying formula (\ref{fpnk}) we get the equality (\ref{fhnk}) when $k<n$.

Let $k=l=n$. If $K_\xi \in f (\langle n,n\rangle )$, then $|K_\xi |=n^2$ and since $\langle n,n\rangle $ is the maximal element in $V_n$, then

$${\rm h} ( n,n )=\frac{\displaystyle |\mathcal{S}_n |-\sum_{\langle r,t \rangle \ne \langle n,n\rangle } rn\, {\rm h} (\langle s,t\rangle  )}{n^2} = $$

$$=\frac{1}{n} \left( (n-1)! \left( \frac{n}{n} \right)^{n-1}
-\sum_{\langle r,t\rangle \prec \langle n,n\rangle } r\, {\rm h}
(\langle r,t\rangle ) \right)$$

\hfill $\Box$

\begin{corollary} \cite{ki} \label{prime}
When $n$ is a prime number
\begin{equation} \label{Nprime}
|{\mathcal{S}_n }_{/_{\stackrel{\sigma}{\sim}}} |= \frac{(n-1)! +(n-1)^2}{n}
\end{equation}
\end{corollary}

Proof. When $n$ is a prime number according to lemma \ref{lm3},  $V_n$ consists only of  the vertices $\langle 1,l\rangle$ and $\langle n,n\rangle$, where $l\in \{ 1,2,\ldots ,(n-1)\}$. From the corollary \ref{kk=1} it follows that ${\rm h}(n,1)=1$, and when $k=n$ according to (\ref{fhnk}) there is
$${\rm h} ( n,n ) =\frac{(n-1)! -(n-1)}{n}.$$

Then we get:
$$|{\mathcal{S}_n }_{/_{\stackrel{\sigma}{\sim}}} |= (n-1)\cdot 1 + \frac{(n-1)! -(n-1)}{n} =$$
$$\frac{(n-1)! +(n-1)^2}{n} .$$

\hfill $\Box$

Since the number of the equivalence classes according to the given equivalence relation is whole number, then from the corollary \ref{prime} follows the proof of the famous in the number theory Wilson's theorem:

\begin{corollary}  \label{wilson} {\bf (Wilson's Theorem)}
If $n$ is a prime number, then
\begin{equation}
(n-1)! +1 \equiv 0 \pmod{n}
\end{equation}
\hfill $\Box$
\end{corollary}

\begin{theorem}\label{varphink}
The number of the vertices of kind $\langle k,l\rangle \in V_n$ is equal to  $\displaystyle \varphi \left( \frac{n}{k}\right)$, where $\varphi (m)$  is the Euler totient function.
\end{theorem}

Proof. At $k=1$ and $k=n$, the assertion is obvious. Let $1<k<n$. We examine the sets:
\begin{equation}\label{W(n,k)}
W(n,k)=\left\{ \langle k,l\rangle \in V_n \right\} \subset V_n
\end{equation}

\begin{equation}\label{Phi(m)}
\Phi (m)=\left\{ s\in \mathbb{N} \; \left| \; 1\le s<m,\; \mathrm{GCD}(s,m)=1 \right. \right\}
\end{equation}

Obviously $\left| \Phi (m)\right| =\varphi (m)$.

Firstly we demonstrate that for every $\displaystyle t\in \Phi \left( \frac{n}{k} \right) $, i.e. $\displaystyle 1\le t <\frac{n}{k}$ and $\displaystyle \mathrm{GCD} \left( t,\frac{n}{k} \right) =1$ the vertex $\langle k, kt \rangle \in W(n,k)\subset V_n$ exists.

 i. Let  $\mathrm{GCD} (t,n)=1$. Then the vertex  $\langle 1,t\rangle \in V_n$ exists. As $\displaystyle t<\frac{n}{k}$, then $kt<n$, from where it follows that the  vertex $\langle k,kt\rangle \in W(n,k)$ exists.

 ii. Let  $\mathrm{GCD} (t,n)=q\ge 2$. Then  $n$ can be presented in kind  $n=q_1 n_1$, where  $q_1 \ge q\ge2$, $q|q_1$, so that if  $p$ is a prime divisor of $n$, then  either $p|n_1$, but $p\nmid q \Rightarrow p\nmid t$, or $p\nmid n_1$, but $p| q \Rightarrow p|t$.

 We assume $s=n_1 +t$. Obviously $\mathrm{GCD} (s,n)=1$. Then
 $$\displaystyle s=n_1 +t<n_1 +\frac{n}{k} = \frac{n}{q_1} +\frac{n}{k} = n\left( \frac{1}{q_1} +\frac{1}{k} \right) \le n\left( \frac{1}{2} +\frac{1}{2} \right) =n.$$

 Consequently $s<n$, i.e. $s\in \Phi (n)$, from where it follows that $\langle 1,s\rangle \in V_n$.

As  $\displaystyle \mathrm{GCD} \left( t,\frac{n}{k} \right) =1$, then it is easily to see that $q_1 |k$, i.e. a positive integer $k_1$ exists, such that $k=q_1 k_1$. Then:
$$ks=k\left( n_1 +t\right) = q_1 k_1 \frac{n}{q_1} +kt\equiv kt\; (\mathrm{mod}\; n)$$

Moreover   $\displaystyle kt< k\frac{n}{k} =n$. Consequently   $\langle k, kt\rangle \in W(n,k)$ in this case too.

In this way we demonstrate that a mapping $\displaystyle \rho : \Phi \left( \frac{n}{k} \right) \to W(n,k)$ from  $\displaystyle  \Phi \left( \frac{n}{k} \right) $ into $ W(n,k)$ exists such that for every $\displaystyle t\in \Phi \left( \frac{n}{k} \right)$ the vertex  $\langle k, kt\rangle \in W(n,k)$ to correspond.

Let  $\displaystyle t_1 ,t_2 \in \Phi \left( \frac{n}{k} \right)$, $t_1 \ne t_2$. Then   $\displaystyle 1\le t_1 <\frac{n}{k}$  and  $\displaystyle 1\le t_2 <\frac{n}{k}$, from where $1<k\le kt_1  <n$ and $1<k\le kt_2  <n$. Consequently   $kt_1 \not\equiv kt_2 \; (\mathrm{mod}\; n)$. Consequently the mapping  $\displaystyle \rho : \Phi \left( \frac{n}{k} \right) \to W(n,k)$ is injection.

Let  $\langle k,l\rangle \in W(n,k)$. It follows from lemma \ref{lm3} that a positive integer $l_1$ exists such that  $l=kl_1$.  As  $kl_1 <n$, then $\displaystyle l_1 < \frac{n}{k}$We suppose that $\displaystyle l_1 \not\in \Phi \left( \frac{n}{k} \right)$, i.e. we suppose that a integer  $p>1$ exists such that $p|l_1$ and $\displaystyle p| \left( \frac{n}{k} \right) $. Consequently integers  $u\ge 1$ and $v\ge 1$ exist such that $\displaystyle \frac{n}{k} =pu \Rightarrow n=kpu$ and $l_1 = pv$. Then   $ul=ukl_1 =ukpv =nv \equiv n\; (\mathrm{mod} \; n)$. Then it  follow from  $\langle k,l\rangle \in V_n$ that  $\langle ku,n\rangle \in V_n$, which is impossible as   $ku<n$ and the only vertex of kind    $s,n$ is at $s=n$. Consequently the mapping   $\displaystyle \rho : \Phi \left( \frac{n}{k} \right) \to W(n,k)$ is surjection.

Consequently a bijective mapping from    $\displaystyle  \Phi \left( \frac{n}{k} \right) $ into $ W(n,k)$ exists, from where it follows that
$$\left| W(n,k)\right| =\left| \Phi \left( \frac{n}{k} \right) \right| =\varphi \left( \frac{n}{k} \right).$$

\hfill $\Box$

\begin{corollary}  \label{wwwwww}
Let $n\in \mathbb{N}$. Then the number of all vertices of the graph  $\Gamma_n =\left( V_n , \Psi\right)$ is equal to  $n$.
\end{corollary}

The proof follows from well known equality
$$\displaystyle \sum_{d|n} \varphi (d) =n.$$
\hfill $\Box$

\begin{corollary}  \label{xxxson}

\begin{equation}\label{QnQnQ}
\left| Q_n \right| = \sum_{k|n} h(n,k) \varphi \left( \frac{n}{k} \right)
\end{equation}
\hfill $\Box$
\end{corollary}

\section{Matrix representation}

The formulas (\ref{hnk}) and (\ref{QnQnQ})  give us an effective algorithm for the manual calculation of  $|Q_n |$. Construction of the graph  $\Gamma_n$ is necessary for this object.  This approach gives relatively good results at relatively small values of $n$  as the experience of author has been shown. With increase of  $n$ probability of errors increases repeatedly, because of the number of classes increases exponentially, according to formulas (\ref{hnk}) and (\ref{QnQnQ}).

In this section we will describe an algorithm for solving the combinatorial problem for finding the cardinality of $Q_n$ and which is suitable for computer implementation. This algorithm is based on proven above statements. For this purpose we transform the formula \ref{fhnk} in the following recursive form

\begin{equation}\label{hnktau}
{\rm h}(n,k) =\frac{1}{k} \left( (k-1)! \left( \frac{n}{k} \right)^{k-1} - \sum_{r|k,\; r<k }  r\tau (n,k,r) {\rm h}(n,r) \right) \quad \mathrm{for} \; k>1
\end{equation}
and
\begin{equation}\label{hnktau1}
{\rm h}(n,1) =1
\end{equation}
where the function $\tau(n,k,r)$  gives the number of vertices  $\langle r,s\rangle \in V_n$, such that  $\langle r,s\rangle \prec \langle k,k\rangle \in V_n$.

For using the formula (\ref{QnQnQ}) as we avail ourselves of formulas (\ref{hnktau}) and (\ref{hnktau1}) it is convenient to fill the elements of the matrix  $M=\left( m_{ij} \right)$, consisting of 4 rows and $q$  columns, where $q$  is equal to the number of all positive divisors of the integer $n$  (including 1 and  $n$).

In the first row we write down all integers dividing without remainder the parameter $n$.

We fill the second row of $M$  with the values of the Euler function  $\displaystyle \varphi \left( \frac{n}{k} \right)$, where $k$  has been taken from the corresponding component of the first row of  $M$. The function $\varphi (m)$  can be realized as we use a variety of the algorithm, known as "the Sieve of Eratosthenes" \cite{mirchevNumber,reingold}.

We fill the third row of  $M$ with the values of the function $h(n,k)$, according to formulas (\ref{hnktau}) and (\ref{hnktau1}), where  $k$ has been taken from the corresponding component of the first row of $M$. For the calculation of the function  $h(n,k)$, must first be calculated the function  $\tau (n,k,r)$ which gives the number of vertices  $\langle r,s\rangle \in V_n$, such that  $\langle r,s\rangle \prec \langle k,k\rangle \in V_n$.

We receive the fourth row of $M$   as multiply component by component the second and the third rows.

\begin{table}
\begin{center}
\caption{Matrix $M$ at $n=12$.}\label{atn=12}

\begin{tabular}{|c||c|c|c|c|c|c|} \hline
$k|n$ & 1 & 2 & 3 & 4 & 6 & 12 \\ \hline
$\displaystyle \varphi \left( \frac{n}{k} \right)$ & 4 & 2 & 2 & 2 & 1 & 1\\ \hline
$h(n,k)$ & 1 & 2 & 10 & 39 & 628 & 3326054\\ \hline
$m_{2j}*m_{3j}$ & 4 & 4 & 20 & 78 & 628 & 3326054 \\
\hline
\end{tabular}
\end{center}

\end{table}

We receive the final result for the number $|Q_n |$  of the equivalence classes of $\mathcal{S}_n$  by the considered equivalence relation as we add the elements of the last row of the matrix  $M$  according to formula (\ref{QnQnQ})

For example at $n=12$  the filled matrix $M$  is shown in table 1 and so
$$|Q_{12} |=4+4+20+78+628+3326054=3326788 .$$

\begin{table}
\begin{center}
\caption{The number of the equivalence classes in  $\mathcal{S}_n$  at  $2\le n\le 19$.}\label{table2}

\begin{tabular}{|c||c|c|c|c|c|c|c|c|c|} \hline
$n$      & 2 & 3 & 4 & 5 & 6  & 7   & 8   & 9        \\ \hline
$|Q_n |$ & 1 & 2 & 3 & 8 & 24 & 108 & 640 & 4 492 \\
\hline
\end{tabular}

\bigskip
\begin{tabular}{|c||c|c|c|c|c|} \hline
$n$        & 10     &  11     & 12         & 13    \\ \hline
$|Q_n |$   & 36 336 & 329 900 & 3 326 788   & 36 846 288\\
\hline
\end{tabular}

\bigskip
\begin{tabular}{|c||c|c|c|c|} \hline
$n$         & 14          & 15            & 16               \\
 \hline
$|Q_n |$ & 444 790 512 & 5 811 886 656 & 81 729 688 428    \\
\hline
\end{tabular}

\bigskip
\begin{tabular}{|c||c|c|c|c|c|c|c|} \hline
$n$      & 17                & 18                 & 19                   \\  \hline
$|Q_n |$ & 1 230 752 346 368 & 19 760 413 251 956 & 336 967 037 143 596  \\
\hline
\end{tabular}
\end{center}

\end{table}

\section{Conclusion}
We show in table \ref{table2} the values of $|Q_n |$   at $2\le n\le 19$  calculated by a computer program based on the above described algorithm. The software is implemented in C++ programming language, and it is described in \cite{Lilly}.

\end{document}